\documentclass{article}
\usepackage[dvipdfmx]{graphicx}
\usepackage{spconf}
\usepackage{latexsym}
\usepackage{amsmath}
\usepackage{amssymb}
\usepackage{amsfonts}
\usepackage{amsthm}
\usepackage{setspace}
\usepackage[ruled,vlined]{algorithm2e}
\usepackage{url}
\usepackage{cite}

\makeatother

\theoremstyle{plain}

\newtheorem{remk}{Remark}

\newcommand{\argmin}{\mathop{\rm argmin}\limits}

\graphicspath{{./images/}}
\def\x{{\mathbf x}}

\def\y{{\mathbf y}}

\def\h{{\mathbf h}}
\def\f{{\mathbf f}}

\def\A{{\mathbf A}}
\def\B{{\mathbf B}}

\def\F{{\mathbf F}}
\def\G{{\mathbf G}}

\def\R{{\mathbb R}}
\def\N{{\mathbb N}}
\def\I{{\mathbf I}}

\def\W{{\mathbf W}}

\def\X{{\mathbf X}}
\def\Y{{\mathbf Y}}

\def\Lmath{{\mathcal L}}

\def\Mmath{{\mathcal M}}

\def\Xbcal{{\boldsymbol{\mathcal X}}}
\def\Ybcal{{\boldsymbol{\mathcal Y}}}

\def\Ebcal{{\boldsymbol{\mathcal E}}}

\DeclareMathOperator{\prox}{prox}

\title{Efficient Constrained Tensor Factorization by\\Alternating Optimization with Primal-Dual Splitting}

\name{Shunsuke Ono${}^{\dag}$ and Takuma Kasai${}^{\dag\dag}$
\thanks{%We thank the reviewers for their careful reading and valuable comments.
The work was partially supported by JSPS Grants-in-Aid (17K12710) and JST-PRESTO.}}
\address{%Tokyo Institute of Technology%, Imaging Science $\&$ Engineering Laboratory, Kanagawa, Japan\\
$^{\dag}$Tokyo Institute of Technology \;\;\;\;$^{\dag\dag}$RIKEN
}

\begin{document}
\maketitle

\begin{abstract}
Tensor factorization with hard and/or soft constraints has played an important role in signal processing and data analysis.
However, existing algorithms for constrained tensor factorization have two drawbacks: (i) they require matrix-inversion; and (ii) they cannot (or at least is very difficult to) handle \textit{structured} regularizations.
We propose a new tensor factorization algorithm that circumvents these drawbacks.
The proposed method is built upon alternating optimization, and each subproblem is solved by a primal-dual splitting algorithm, yielding an efficient and flexible algorithmic framework to constrained tensor factorization.
The advantages of the proposed method over a state-of-the-art constrained tensor factorization algorithm, called AO-ADMM, are demonstrated on regularized nonnegative tensor factorization.
\end{abstract}

\section{Introduction}\label{intro}
Tensor factorization techniques have been extensively studied and applied not only to signal processing and machine learning problems, including signal analysis and blind source separation \cite{cichocki2009nonnegative}, dimensionality reduction and learning latent variable models \cite{symeonidis2008tag,anandkumar2014tensor}, but also to scientific problems in chemometrics \cite{smilde2005multi} and neuroscience \cite{morup2006parallel}.
Recent comprehensive reviews on tensor factorization can be found in \cite{cichocki2015tensor,sidiropoulos2017tensor}
%However, despite the fact that low-rank tensor factorization is essentially unique under mild conditions \cite{sidiropoulos2000uniqueness}, determining tensor rank is NP-hard, and the best low-rank approximation of a higher rank tensor may not even exist.

In the so-called \textit{canonical polyadic decomposition} (CPD) model, also known as the \textit{parallel factor analysis} (PARAFAC) model \cite{carroll1970analysis,PARAFAC}, a tensor is decomposed into a sum of the lowest possible number of rank-1 tensors, where a rank-one tensor consists of an outer product of vectors.
Although CPD is essentially unique under relatively mild conditions \cite{sidiropoulos2000uniqueness},
hard and/or soft constraints on factors, such as nonnegativity, sparsity, smoothness and so on, are very useful for restoring identifiability,
improving estimation accuracy, %in relatively challenging cases (e.g., low-SNR),
ensuring interpretability of the results, %(e.g., nonnegative values),
and fixing ill-posedness \cite{sidiropoulos2017tensor}.
On the other hand, although there exist a bunch of tensor factorization algorithms, most of them are designed for unconstrained cases or customized for handling a specific constraint, as pointed out in \cite{huang2016flexible}.

Very recently, a tensor factorization algorithm
that can easily and naturally incorporate various types of constraints has been proposed \cite{huang2016flexible}.
This method is based on the \textit{alternating optimization}, i.e., updating each variable (factor matrix) by solving the corresponding subproblem in a cyclic fashion,
which is the standard approach adopted in many other tensor factorization algorithms.
What is different is how to solve each subproblem: the method adopts the \textit{alternating direction method of multipliers} (ADMM) \cite{ADMM1,DRS2,ADMMBoyd}.
ADMM can efficiently solve nonsmooth convex optimization problems with the help of proximal splitting techniques \cite{TechCombettes}.
%What is different is how to solve each subproblem: the one proposed in \cite{} adopts the proximal gradient algorithm \cite{}, wheares the other one \cite{} uses the alternating direction method of multipliers \cite{}.
Thereby, this method, which the authors named AO-ADMM, deals with involved subproblems that have no closed-form solutions.

However, there are two things to be improved.
One is that AO-ADMM requires matrix-inversion at each iteration of ADMM.
The authors of \cite{huang2016flexible} suggest to alleviate this computational difficulty by Cholesky decomposition and back-substitution, but it is still a bottleneck.
The other is that a class of \textit{structured} regularization, i.e., the composition of a simple regularization function and a linear operator, cannot be (or at least is very difficult to be) used as a soft constraint.
Representative examples include the overlapping group lasso \cite{OverlapGroupLasso} and the total variation \cite{ROF},
which would be useful in many applications, as the lasso and the quadratic variation having been commonly used \cite{allen2012sparse,liu2012sparse,yokota2016smooth}.
%Representative examples include the overlapping group lasso \cite{OverlapGroupLasso}, which effectively promotes group sparsity, and the total variation \cite{ROF}, which imposes piece-wise smoothness.
Note that these drawbacks are common to other constrained tensor factorization algorithms.

To circumvent these drawbacks, we propose a new algorithmic framework for constrained tensor factorization.
The proposed method is built upon alternating optimization as well as AO-ADMM, but the essential difference is that each subproblem is solved by a \textit{primal-dual splitting algorithm} \cite{PDCondat,PDSVu}.
It can solve nonsmooth convex optimization problems involving linear operators without inversion and has been applied to signal and image processing problems, e.g., \cite{PDCondatletter,CVPR2014,PDAndernonconvex,PDSICIP,TSP2015Hierarchical,SPL2017}.
Thus, incorporating it into alternating optimization yields a more efficient and flexible algorithm for constrained tensor factorization than AO-ADMM.
The advantages of the proposed method are demonstrated on regularized nonnegative tensor factorization, where aside from its efficiency, we empirically show that our method achieves better factorization in terms of mean squared error.

\section{Preliminaries}
\subsection{Canonical Polyadic Decomposition (CPD)}
For brevity, we focus on third-order tensors in this paper but everything naturally generalizes to higher-order tensors.
We denote a tensor by bold calligraphic letters like $\Xbcal\in\R^{N_1\times N_2 \times N_3}$,
a matrix by bold capital letters like $\X\in\R^{N_1\times N_2}$,
and a vector by bold letters like $\x\in\R^N$.

We denote %the Kronecker product of $\X\in\R^{N_1\times N_2}$ and $\Y\in\R^{M_1\times M_2}$ by $\X\otimes\Y\in\R^{M_1N_1\times M_2N_2}$,
the Khatri-Rao product (column-wise Kronecker product) of $\X$ and $\Y$ with the same number of columns (i.e.. $M_2=N_2$) by $\X\odot\Y\in\R^{M_1N_1\times N_2}$,
and the tensor product (outer product) of $\x\in\R^{N}$ and $\y\in\R^{M}$ by $\x\circledcirc\y\in\R^{N\times M}$.
Note that the tensor product of three vectors yields a rank-1 third-order tensor.
The mode-$d$ matricization of $\Xbcal$ is a matrix of size $\prod N_{i\neq d} \times N_d$, denoted by $\X_{(d)}$.
Each row of $\X_{(d)}$ is a vector obtained by fixing all the indices of $\Xbcal$ except the $d$-th one, and the matrix is
formed by stacking these row vectors by traversing the rest of the indices from 3 back to 1.

CPD of $\Xbcal$ is given by
\begin{equation}
\Xbcal = \sum_{r=1}^R \circledcirc_{d=1}^3\f_d^{(r)},
\end{equation}
where $\f_d^{(r)}\in\R^{N_d}$ ($d=1,2,3$, $r=1,\ldots,R$) are the factors of $\Xbcal$, and $R>0$ is the rank of $\Xbcal$ that is the minimum number of rank-1 tensors required to represent $\Xbcal$ as their sum.
Note that a predetermined $0<\tilde R\leq\min_d\{N_d\}$ is used in practice instead of the true tensor rank $R$ since computing $R$ in CPD is NP-hard \cite{hillar2013most}.
We denote the above relation by
\begin{equation}\label{CPDmat}
\Xbcal = [\F_d]_{d=1}^3,
\end{equation} where $\F_d:=(\f_d^{(1)}\cdots\f_d^{(R)})\in\R^{N_d\times R}$ is the factor matrix of the $d$-th mode.
With the help of this notation, we can express CPD in a matricized form:
\begin{equation}\label{CPDmatricize}
\X_{(d)} = (\odot_{i\neq d}\F_i)\F_d^\top.
\end{equation}

\subsection{Primal-Dual Splitting Algorithm}
Let $\Gamma_0(\R^N)$ be the set of all proper lower semicontinuous convex functions on $\R^N$.
Consider convex optimization problems of the form:
\begin{equation}\label{prob:PDS}
\min_{\x\in\R^N} f(\x) + g(\x) + h(\Lmath\x),
\end{equation}
where $f,g\in\Gamma_0(\R^N)$ ($f$ is $\beta$-Lipschitz differentiable with $\beta>0$), $h\in\Gamma_0(\R^K)$, and $\Lmath:\R^{N}\rightarrow\R^{K}$ is a linear operator.
Also, let us introduce the notion of the \textit{proximity operator} of index $\gamma>0$ of $f\in\Gamma_0(\R^N)$
as follows:
\begin{equation}\label{prox}
\prox_{\gamma f}:\R^N\rightarrow\R^N:\x\mapsto\argmin_{\y}f(\y)+\frac{1}{2\gamma}\|\y-\x\|^2.\hspace{-1mm}
\end{equation}

A primal-dual splitting (PDS) algorithm \cite{PDCondat} solves Prob.~\eqref{prob:PDS} by
the following iterative procedure: for any $\x^{(0)}\in\R^N$, $\y^{(0)}\in\R^K$, and $\gamma_1,\gamma_2>0$ satisfying $\gamma_1(\frac{\beta}{2}+\gamma_2\|\Lmath^*\Lmath\|)< 1$ ($\Lmath^*$ is the adjoint operator of $\Lmath$, and $\|\cdot\|$ denotes the operator norm), iterate
\begin{equation}\label{PDS}
\hspace{-0.5mm}\left\lfloor\begin{array}{l}
\x^{(n+1)}:=\prox_{\gamma_1 g}(\x^{(n)}-\gamma_1(\nabla f(\x^{(n)}) + \Lmath^*(\y^{(n)}))),\\
\y^{(n+1)}:=\prox_{\gamma_2 h^*}(\y^{(n)} + \gamma_2\Lmath(2\x^{(n+1)}-\x^{(n)})),
\end{array}\right.\hspace{-0.5mm}
\end{equation}
where $h^*$ is the convex conjugate function of $h$,
and $\gamma_1$ and $\gamma_2$ can be seen as the stepsizes.
We note that the proximity operator of $h^*$ is available via that of $h$ as
\begin{equation}\label{dualprox}
\prox_{\gamma h^*}(\x)=\x - \gamma\prox_{\gamma^{-1}h}(\gamma^{-1}\x)
\end{equation}
(see, e.g., \cite[Theorem~14.3(ii)]{Combettesbook}). %, and $\gamma_1, \gamma_2>0$ satisfy $\gamma_1\gamma_2\|A\|_{op}^2\leq 1$
Under some mild condition on $g$, $h$, and $\Lmath$,
the sequence $(\x^{(n)})_{k\in\N}$ converges to an optimal solution of Prob.~\eqref{prob:PDS}.

\section{Proposed Method}
\subsection{Problem Formulation}
Consider the following data observation model:
\begin{equation}\label{model}
\Ybcal = \Mmath(\Xbcal + \Ebcal),
\end{equation}
where $\Ybcal$ is an observed data stored as a tensor possibly with missing data, $\Xbcal$ is a true tensor data, $\Mmath$ is a self-adjoint idempotent linear operator that specifies the observed entries in $\Ybcal$, i.e., zeroing out the entries whose indices correspond to missing data, and $\Ebcal$ is an additive noise.

Then, wth the notation in \eqref{CPDmat}, we formulate constrained tensor factorization as a generic optimization problem:
\begin{align}
\min_{\F_1,\F_2,\F_3} \frac{1}{2}\|\Ybcal - \Mmath([\F_d]_{d=1}^3)\|_F^2 + \sum_{d=1}^3 h_d(\Lmath_d(\F_d^\top))\nonumber\\
\mbox{ s.t. } \F_d^\top\in C_d \;\; (d=1,2,3),\label{mainprob}
\end{align}
where $\|\cdot\|_F$ is the Frobenius norm of a tensor, $h_d\circ \Lmath_d$ is a regularization function (soft constraint) for the $d$-th factor matrix consisting of a (possibly nonsmooth) convex function $h_d$ and a linear operator $\Lmath_d$, and $C_d$ is a closed convex set representing a hard constraint on the $d$-th factor matrix.
Here we assume that the proximity operator
%\footnote{The \textit{proximity operator} of a convex function $f$ is defined by
%	\begin{equation*}\label{prox}
%	\prox_{\gamma f}:\R^N\rightarrow\R^N:\x\mapsto\argmin_{\y}f(\y)+\tfrac{1}{2\gamma}\|\y-\x\|^2.
%	\end{equation*}}
of $h_d$ and the metric projection\footnote{The \textit{metric projection} onto a closed convex set $C$ is given by
	\begin{equation*}\label{proj}
	P_C:\R^N\rightarrow\R^N:\x\mapsto\argmin_{\y\in C}\|\y-\x\|^2.
	\end{equation*}} onto $C_d$ are efficiently computable.

This formulation covers various existing constrained tensor factorization problems.
A typical example would be nonnegative tensor factorization \cite{bro1997fast,shashua2005non},
which is recovered by setting $h_d:=0$ and $C_d:=\R_+^{N_d\times R}$ ($\R_+$ denotes the set of all nonnegative real numbers).
Another example is $\ell_1$-regularized tensor factorization \cite{allen2012sparse,liu2012sparse} promoting the sparsity of factors,
which corresponds to penalizing factors by the $\ell_1$ norm, i.e., $h_d:=\|\cdot\|_1$.
These constraints have shown to be very useful for the reasons mentioned in Sec.~\ref{intro}.

\subsection{Optimization}
Since Prob.~\eqref{mainprob} is nonconvex due to the multi-linearity of CPD,
it is very difficult to solve the problem directly.
As a remedy, the alternating optimization is commonly used: each factor matrix $\F_d$ is updated in a cyclic fashion.
Specifically, for each $\F_d$, we solve the following subproblem:
\begin{align}
\min_{\F_d} \frac{1}{2}\|\Y_{(d)} - \Mmath((\odot_{i\neq d}\widetilde{\F}_i)\F_d^\top)\|_F^2 + h_d(\Lmath_d(\F_d^\top))\nonumber\\
\mbox{ s.t. } \F_d^\top\in C_d,\label{subprob}
\end{align}
where we use the matricized form in \eqref{CPDmatricize}, and $\widetilde{\F}_i$ ($i\neq d$) are the fixed factor matrices except the mode $d$.

Now we can see that Prob~\eqref{subprob} is convex but is still a tough problem because of the nonsmoothness.
Thus, we propose to approximately solve the problem by few iterations of the primal-dual splittihg algorithm in \eqref{PDS}.
To this end, first, we introduce the indicator function of $C_d$, defined by $\iota_{C_d}(\x):=0$, if $\x\in C$; $\iota_{C_d}(\x):=\infty$, otherwise.
It should be noted that the proximity operator of the indicator function of $C_d$ is equivalent to the metric projection onto $C_d$.
Second, by letting $\W:=\odot_{i\neq d}\widetilde{\F}_i$ and $\F:=\F_d^\top$, Prob~\eqref{subprob} can be rewritten as
\begin{align*}
\min_{\F} \frac{1}{2}\|\Y_{(d)} - \Mmath(\W\F)\|_F^2 + \iota_{C_d}(\F) + h_d(\Lmath_d(\F)).
\end{align*}

Let us define
\begin{align*}
f(\F)&:=\frac{1}{2}\|\Y_{(d)} - \Mmath(\W\F)\|_F^2,\\
g(\F)&:=\iota_{C_d}(\F),\\
h(\G)&:=h_d(\G) \mbox{ and } \Lmath:=\Lmath_d.
\end{align*}
Since the squared loss term is $\beta$-Lipschitz differentiable with
\begin{align*}
\beta&=\|\W^\top\Mmath^*(\Mmath(\W))\|=\|\W^\top\Mmath(\W)\| \\
&\leq\|\W^\top\|\|\Mmath\|\|\W\|=\|\W^\top\W\| \;\;(\because \|\Mmath\|=1),
\end{align*}
and the proximity operators of $\iota_{C_d}$ and $h_d$ are available from the assumptions on $C_d$ and $h_d$, we can derive an iterative algorithm for solving Prob.~\eqref{subprob} based on the primal-dual splitting algorithm in \eqref{PDS} as follows: for any $\F^{(0)}$, $\G^{(0)}$, and $\gamma_1,\gamma_2>0$ satisfying
\begin{equation}
\gamma_1\left(\frac{\|\W^\top\W\|}{2}+\gamma_2\|\Lmath_d^*\Lmath_d\|\right)< 1,\label{pdsineq}
\end{equation}
set $\A:=\W^\top\Mmath(\W)$ and $\B:=\W^\top\Y_{(d)}$, and iterate
\begin{equation}\label{PDSsub}
\hspace{-0.5mm}\left\lfloor\begin{array}{l}
\F^{(n+1)} := P_{C_d}(\F^{(n)}-\gamma_1(\A\F^{(n)} - \B + \Lmath_d^*(\G^{(n)}))),\\
\G^{(n)} \leftarrow \G^{(n)} + \gamma_2\Lmath_d(2\F^{(n+1)}-\F^{(n)}),\\
\G^{(n+1)} := \G^{(n)} - \gamma_2\prox_{\frac{1}{\gamma_2} h_d}(\frac{\G^{(n)}}{\gamma_2}).
\end{array}\right.\hspace{-0.5mm}
\end{equation}

\begin{remk}[Comparison with AO-ADMM \cite{huang2016flexible}]\normalfont
Clearly, there is no matrix inversion in \eqref{PDSsub}.
Specifically, whereas AO-ADMM requires the inversion of $\A + \I$ (see Algorithm~1 in \cite{huang2016flexible}), our algorithm only needs to compute the multiplication of $\A$ and $\F^{(n)}$.
The authors of \cite{huang2016flexible} suggest to use Cholesky decomposition and back-substitution for efficiently computing this inversion, but we will see in the next section that our algorithm outperforms AO-ADMM in terms of CPU time.
We also remark that at the update of $\G$, our algorithm just computes the proximity operator of $h_d$, i.e., the linear operator $\Lmath_d$ is decoupled.
This is another big difference from AO-ADMM: it requires to compute the proximity operator of $h_d\circ\Lmath_d$, which does \textit{not} have a closed-form solution in general even if that of $h_d$ does.
Indeed, such a situation arises in the following cases: the overlapping group lasso, i.e., $h_d := \|\cdot\|_1$ and $\Lmath_d$ is an operator that replicates overlapping variables, and the total variation, i.e., $h_d := \|\cdot\|_1$ and $\Lmath_d$ is a discrete difference operator.

It should be noted that AO-ADMM can be used not only for the squared loss function (the first term in \eqref{subprob}) but also other loss functions such as $\ell_1$ norm and Kullback-Leibler divergence.
Although we only describe the squared loss case in this paper, our approach can also handle other cases by letting $f(\F):=0$, $h(\G_1, \G_2) := l(\G_1) +\h_d(\G_2)$ and $\Lmath:=(\Mmath\W, \Lmath_d)$ in \eqref{prob:PDS}, where $l$ is a loss function.
\end{remk}

%\begin{remk}[Convergence of the proposed algorithm]
%\end{remk}

Finally, by incorporating \eqref{PDSsub} into alternating optimization, we obtain a new algorithm for solving the generic constrained tensor factorization problem formulated in \eqref{mainprob}.
The whole algorithm is shown in Algorithm~\ref{alg:AO-PDS}.
\begin{algorithm}[h]
	\LinesNumbered
	\SetKwInOut{Input}{input}
	\SetKwInOut{Output}{output}
	\label{alg:AO-PDS}
	\caption{Proposed method for solving \eqref{mainprob} (AO-PDS)}
	\Input{$\F_d^{(0)}$, $\G_d^{(0)}$, $\Y_{(d)}$ and $\|\Lmath_d^*\Lmath_d\|$ $(d=1,2,3)$}
	\While{A stopping criterion is not satisfied}{
		\For{$d=1$ to $3$}{
			$\W^{(k)}:=\odot_{i\neq d}\F_i^{(k)}$\;
			$\A^{(k)}:=\W^{(k)\top}\Mmath(\W^{(k)})$\;
			$\B^{(k)}:=\W^{(k)\top}\Y_{(d)}$\;
			Compute $\gamma_1$ and $\gamma_2$ by \eqref{gamma1set} and \eqref{gamma2set}\;
			Update $\F_d^{(k)}$ and $\G_d^{(k)}$ using \eqref{PDSsub} initialized with the previous $\F_d^{(k)}$ and $\G_d^{(k)}$\;
			$\F_d^{(k+1)}\leftarrow\F_d^{(k)}$\;
			$\G_d^{(k+1)}\leftarrow\G_d^{(k)}$\;
		}
	$k\leftarrow k+1$\;
	}
	\Output{$\F_d^{(k)}$ $(d=1,2,3)$}
\end{algorithm}

\begin{figure*}[t]
	\begin{center}
		\begin{minipage}{0.33\hsize}
			\includegraphics[width=\hsize]{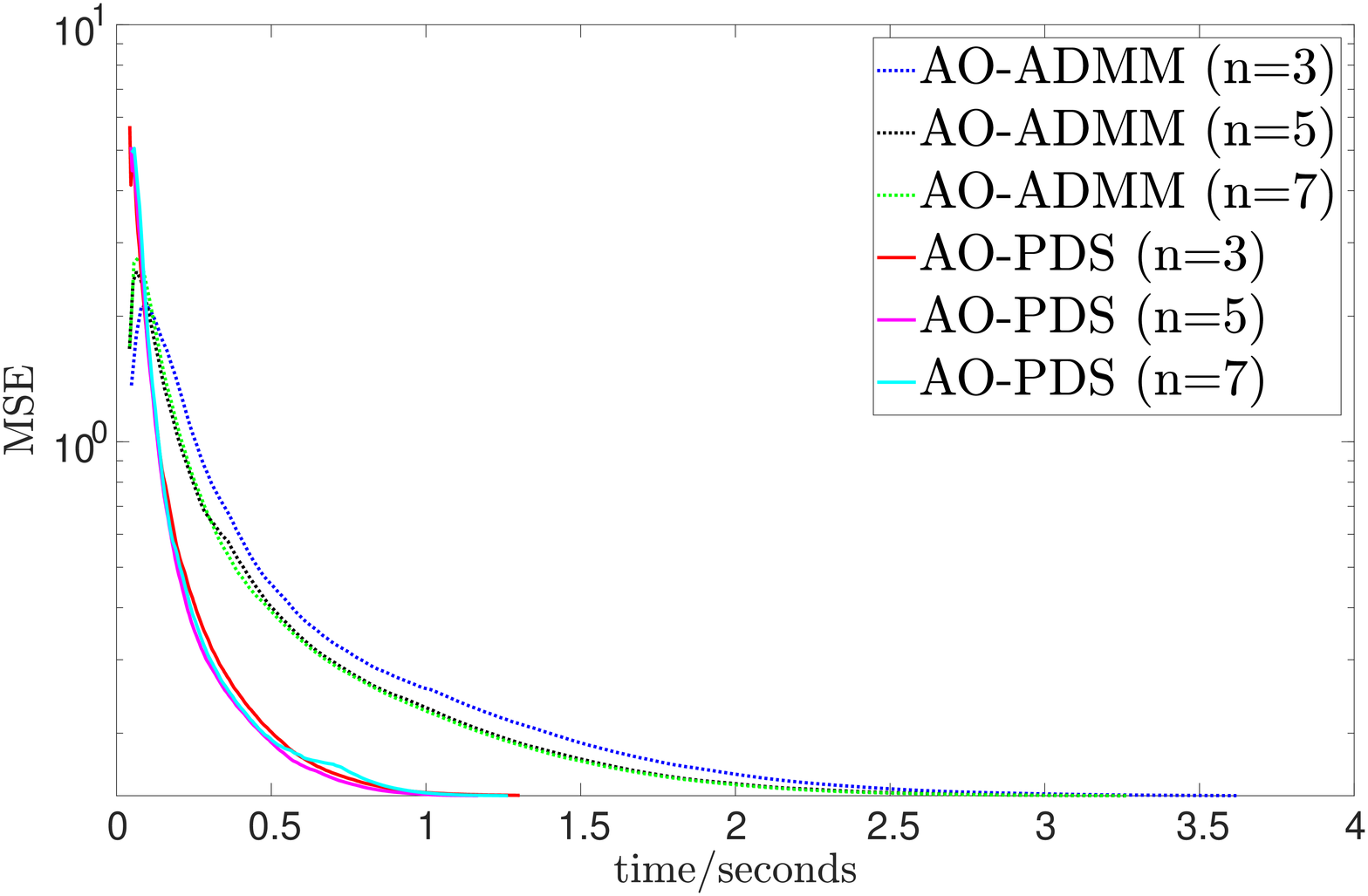}
		\end{minipage}
		\begin{minipage}{0.33\hsize}
	\includegraphics[width=\hsize]{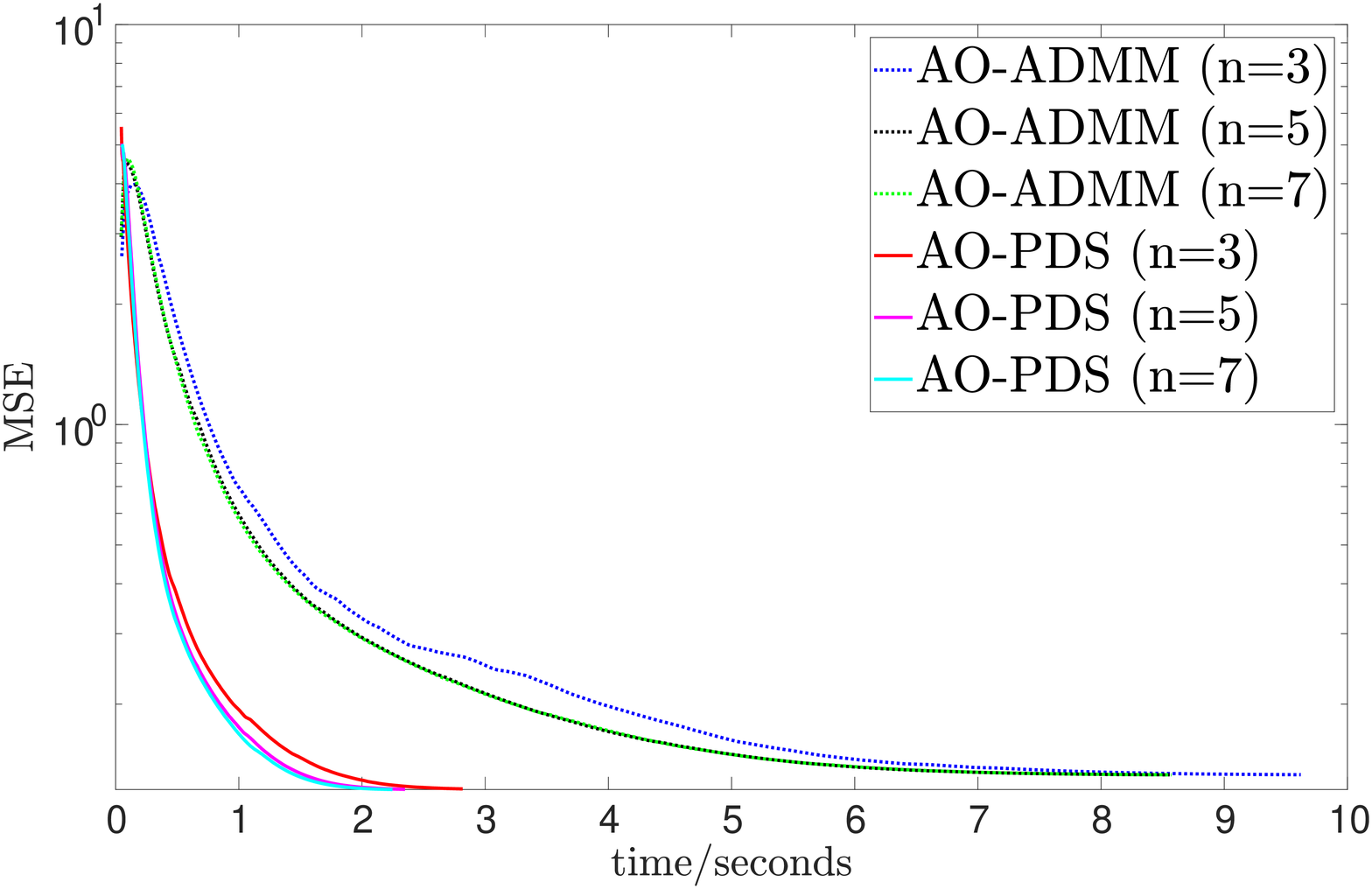}
\end{minipage}
		\begin{minipage}{0.33\hsize}
	\includegraphics[width=\hsize]{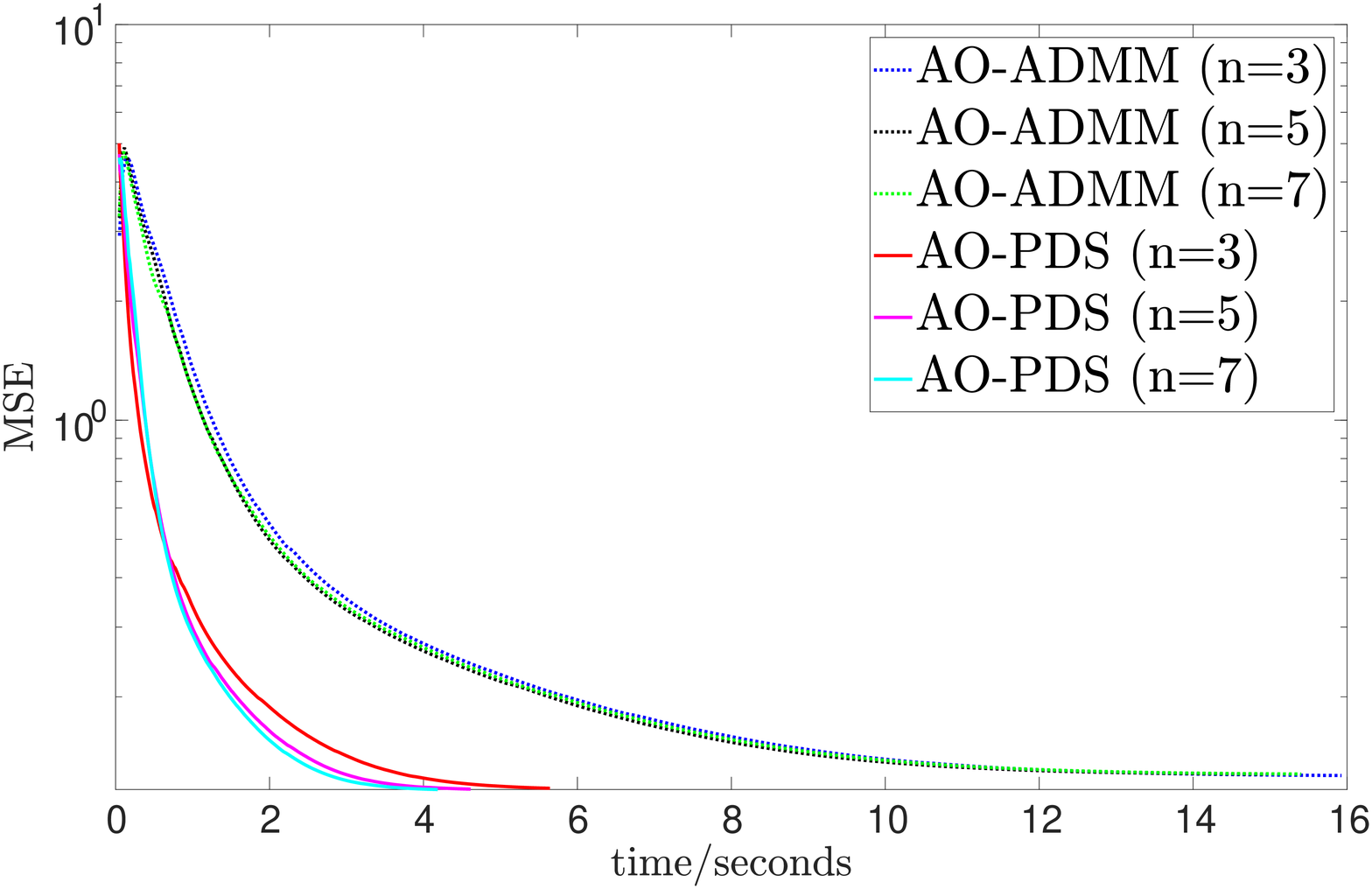}
\end{minipage}

				\vspace{-2mm}
		\caption{Evolution of MSE (logarithmic scale) versus time in seconds on the regularized nonnegative tensor factorization: $R=5$ (left), $R=10$ (center) and $R=15$ (right).}
		\label{fig:converge}
	\end{center}
	\vspace{-7mm}
\end{figure*}

\begin{remk}[Iteration number and stepsizes of \eqref{PDSsub}]\normalfont
	As in the case of AO-ADMM, each factor matrix $\F_d$ and dual variable $\G_d$ are updated in a cyclic fashion in Alg.~\ref{alg:AO-PDS}, so that we can expect that after a number of iterations, these variables obtained in the previous iteration of \eqref{PDSsub} are good initial variables to the current iteration.
	This means that few iterations of \eqref{PDSsub} would be enough for updating $\F_d$ and $\G_d$, which keeps our method computationally efficient.
	Indeed, the numerical experiments in the next section show that a very few number of iterations ($3$ to $7$) are sufficient for empirical convergence with reasonable precision.
	
	In \eqref{PDSsub}, the stepsizes $\gamma_1,\gamma_2>0$ are set to
	\begin{align}
	\gamma_1&:=0.99\frac{2}{\mbox{trace}(\W^\top\W)},\label{gamma1set}\\
	\gamma_2&:=\frac{1}{\gamma_1\|\Lmath_d^*\Lmath_d\|}-\frac{\mbox{trace}(\W^\top\W)}{2\|\Lmath_d^*\Lmath_d\|},\label{gamma2set}
	\end{align}
	respectively.
	This setting satisfies the inequality in \eqref{pdsineq} since the trace of $\W^\top\W$ is an upper bound of $\|\W^\top\W\|$, and it can efficiently be computed at each iteration of Alg.~\ref{alg:AO-PDS}.
	Note that $\|\Lmath_d^*\Lmath_d\|$ can be estimated in advance, and once determined, it can be used for every iteration because $\Lmath_d$ does not change.
\end{remk}

\section{Numerical Experiments}
We examined our method on regularized nonnegative tensor factorization.
Our method was compared with AO-ADMM \cite{huang2016flexible}, where we used the MATLAB code distributed by the authors of \cite{huang2016flexible}.
All experiments were performed using MATLAB (R2017a), on a Windows 8.1 (64bit)
laptop computer with an Intel Core i7 2.6 GHz processor and 16 GB of RAM.

We tested these algorithms on synthetic data.
Specifically, synthetic true tensor data were generated as follows:
for $N_1=N_2=N_3=100$ and $R = 5, 10$ or $15$, the true factor matrices, denoted by $\F_d^{true}$ ($d=1,2,3$), are obtained by
drawing their elements from an i.i.d. uniform distribution on the interval $(0,1)$, and then $80\%$ of the elements of $\F_1^{true}$ are randomly set to $0$, i.e., only $\F_1^{true}$ is sparse.
The observed tensor data $\Ybcal$ is then obtained by \eqref{model}, where the elements of $\Ebcal$ are drawn from an i.i.d. Gaussian distribution with standard deviation $0.1$.
In the experiments, we did not consider missing data, i.e., $\Mmath$ equals to an identity operator.

For soft constraints (regularization), we adopted the $\ell_1$ norm for $\F_1$ and the squared Frobenius norm for $\F_2$ and $\F_3$,
where their hyperparameters were set to $5$ and $2$, respectively.
The nonnnegativity constraint is also imposed on each factor matrix.
Then we measured the evolution of mean squared error (MSE) versus time in seconds on the regularized nonnegative tensor factorization problem solved by AO-ADMM and our algorithm,
where MSE is defined by $\frac{1}{R(N_1+N_2+N_3)}\sum_{d=1}^3 \|\F_d^{true} - \F_d\|_F^2$. %with $\F_d^{true}$ being the true mode-$d$ factor matrix.
Note that the above hyperparameters were hand-optimized in the MSE sense.

For the parameters of AO-ADMM, we used the settings recommended by the authors of \cite{huang2016flexible}.
The stopping criterion of the outer loop of each algorithm is set to $|\mbox{MSE}_k - \mbox{MSE}_{k-1}|<1.0\times 10^{-5}$ ($k$ is the number of iterations of the outer loop).
Note that this criterion can only be used for synthetic data since it uses $\mbox{MSE}$.
For practical situations, one may use the value of the objective function in \eqref{mainprob} as a stopping criterion. 

The results are shown in Figure~\ref{fig:converge},
where $n=3,5,7$ are the number of the inner loop of each algorithm (ADMM or primal-dual splitting).
One can see that our method (AO-PDS) is about \textit{three times faster} than AO-ADMM.
In addition, the best MSE values of AO-PDS are $0.142$, $0.122$ and $0.117$; and those of AO-ADMM are $0.142$, $0.133$ and $0.127$, respectively for $R=5,10,15$, i.e., our method results in \textit{more accurate} factorization in terms of MSE for $R=10,15$.
As expected, we observe that very small number of outer iterations are enough for empirical convergence of AO-PDS.
This would be thanks to the warm-start nature of alternating optimization, as in the case of AO-ADMM.

\section{Conclusion}
We have proposed an efficient and flexible tensor factorization algorithm based on alternating optimization combined with primal-dual splitting.
Our algorithm has two advantages over AO-ADMM, a state-of-the-art tensor factorization algorithm, as follows: (i) it is free from matrix-inversion; and (ii) it can efficiently handle structured regularizations.
Our experimental results on regularized nonnegative tensor factorization not only support our claim that the proposed method is more computationally efficient than AO-ADMM but also revealed that the proposed method achieves better factorization than AO-ADMM in the sense of MSE.

\small{
\bibliographystyle{IEEEbib}
\bibliography{AOPDS.bib}
}

\end{document}